\documentclass[11pt]{amsart}
\usepackage{etex}
\usepackage{diagbox}
\usepackage{array}
\usepackage{graphicx, pictex}
\usepackage{stmaryrd}
\usepackage{changepage}
\usepackage{float}
\restylefloat{table}
\usepackage{multirow}
\usepackage[dvipsnames]{xcolor}
\usepackage{amsmath,amssymb,amsthm,mathrsfs}
\allowdisplaybreaks
\usepackage{yfonts}

\renewcommand{\d}{\partial}

\def\d{\Omega}

\def\du#1#2#3{\overset{#3}{\underset{#2}{#1}}}
\def\Forall{\quad \hbox{ for all }}

\def\M{{\mathcal{M}}}

\newcommand{\tn}[1]{\lVert\kern-1pt\lvert{#1}\rvert\kern-1pt\rVert}
\def\<{{\langle}}
\def\>{{\rangle}}
\def\Forall{\quad \hbox{ for all }}

\def\bb{{\mathbf b}}

\def\d{\Omega}

\def\Forall{\quad \hbox{ for all }}

\def\d{\Omega}

\def\Forall{\quad \hbox{ for all }}

\def\tb#1{{\|\kern-1pt| #1 \|\kern-1pt|}}
\def\nm2#1#2{\|#1\|_{2,\d_{#2}}}

\def\R{\mathbb{R}}

 \theoremstyle{plain}
 \newtheorem{thm}{Theorem}[section]
 \numberwithin{equation}{section} 
 \numberwithin{figure}{section} 
 \theoremstyle{plain}
 
 \theoremstyle{plain}
 \theoremstyle{plain}
 \newtheorem{theorem}[thm]{Theorem}
 \theoremstyle{plain}
 
\theoremstyle{plain}
 \newtheorem{remark}[thm]{Remark}
 \theoremstyle{plain}
 
\def\M{{\mathcal{M}}}

\def\d{{\Omega}}

\def\Forall{\quad \hbox{ for all }}

\def\<{{\langle}}
\def\>{{\rangle}}
\def\R{\mathbb{R}}

\def\bb{{\mathbf b}}

\def\du#1#2#3{\overset{#3}{\underset{#2}{#1}}}

\usepackage{amssymb}

\begin{document}

\title[Explaining oscillatory behavior for CD]
{Explaining oscillatory behavior in  convection-diffusion discretization}

\author{Constantin Bacuta}
\address{University of Delaware,
Mathematical Sciences,
501 Ewing Hall, 19716}
\email{bacuta@udel.edu}






\begin{abstract} For a model convection-diffusion problem, we  address  the presence of oscillatory discrete solutions, and study  difficulties  in recovering standard approximation results for its solution.  We justify  the presence of  non-physical oscillations and propose  ways to eliminate oscillations.  A new approach for error analysis that requires establishing optimal discrete infinity error as a first step is introduced and justified. We emphasize that the  discretization of two dimensional convection dominated problems benefit from the efficient  discretization of the corresponding one dimensional problem along each stream line.
Our results  are useful  in building new and robust  discretizations  for  multi-dimensional convection dominated problems.


\end{abstract}
\maketitle

\section{Introduction}


We consider  the model of a  singularly perturbed  convection diffusion problem: 
Given data $f\in L^2(\Omega)$, find $u \in H^1_0(\Omega)$ such that 
\begin{equation}\label{eq:2d-model}
  \left\{
\begin{array}{rcl}
     -\varepsilon\Delta u +\bb\cdot \nabla u & =\ f & \mbox{in} \ \ \ \Omega,\\
      u & =\ 0 & \mbox{on} \ \partial\Omega,\\ 
\end{array} 
\right. 
\end{equation} 
for a  positive constant  $\varepsilon$ and  a bounded domain  $\Omega\subset \R^d$. We assume  that $\varepsilon \ll1$, and that  the vector  $\bb$ is chosen such that a unique solution exists. 

For the one dimensional case, we further assume that $f$ is  continuous  on $[0, 1]$,  $\bb=1$  and look for  $u=u(x)$  such that
\begin{equation}\label{eq:1d-model}
\begin{cases}-\varepsilon u''(x)+ u'(x)=f(x),& 0<x<1\\ u(0)=0, \ u(1)=0. \end{cases}
\end{equation} 

The model problem \eqref{eq:2d-model}  arises  when solving   heat transfer problems in thin domains, as well as when using small step sizes in implicit time discretizations of parabolic convection diffusion type problems, see \cite{Lin-Stynes12}. The solutions to the model problem \eqref{eq:2d-model}  is characterized by  boundary layers,  see e.g., 
\cite{ EEHJ96,  linssT10, Roos-Schopf15, zienkiewicz2014}. 
Approximating such solutions poses  numerical  challenges due to the $\varepsilon$-dependence of the  stability constants and  of the  error estimates. 
In applications, the model problem \eqref{eq:2d-model} could be  coupled with   a more complex system of equations such as the Navier-Stokes equations through the convection vector  $\bb$. 

The variational formulation of \eqref{eq:2d-model} is: Find $u\in Q:=H_0^1(\Omega)$ such that
\begin{equation} \label{eq:var-cdr}
b(v,u):= (\varepsilon \nabla u, \nabla v)+  (\bb \cdot \nabla u, v)=(f,v) \Forall v\in V:=H_0^1(\Omega).
\end{equation} 
Throughout this paper, $(\cdot, \cdot)$ denotes the scalar or vector $L^2$ inner product.

For the discretization  of \eqref{eq:var-cdr}, we assume that  
$\M_h \subset Q$ and $V_h \subset V$ are  finite element spaces  that are compatible spaces, in the sense that  $b(\cdot, \cdot)$ satisfies  a discrete $\inf-\sup$ condition on $V_h \times \M_h$. 

A general Petrov-Galerkin (PG) discretization of \eqref{eq:var-cdr} with  $dim(V_h) = dim(\M_h)$ is:  Find $u_h\in \M_h $ such that
\begin{equation} \label{eq:var-cdr-h}
b(v_h,u_h)=(f,v_h) \Forall v_h\in V_h.
\end{equation} 

 The Saddle Point Least Square (SPLS) approach  uses an auxiliary variable that represents the residual of the variational formulation \eqref{eq:var-cdr-h} on the test space. This brings in another simple equation involving the residual variable. The method  leads to a square symmetric saddle point system. 
The SPLS discretization of \eqref{eq:var-cdr} is: Find $(w_h, u_h) \in V_h \times \M_h$ such that 
\begin{equation}\label{SPLS4model-h}
\begin{array}{lclll}
a_0(w_h,v_h) & + & b( v_h, u_h) &= (f,v_h ) &\ \Forall  v_h \in V_h,\\
b(w_h,q_h) & & & =0   &\  \Forall  q_h \in \M_h,    
\end{array}
\end{equation}
where $a_0(u, v):= (\nabla u, \nabla v)$. The component $u_h$ of $(w_h, u_h)$ is the {\it SPLS discrete solution} of \eqref{eq:var-cdr}. We can  view the PG discretization as a particular case of SPLS discretization for which $w_h=0$ in 
\eqref{SPLS4model-h}.

The Galerkin method applied to \eqref{eq:var-cdr-h}, with $\M_h=V_h$ being a standard finite element  space of continuous piecewise polynomials on uniform meshes, leads to numerical pollution which translates into  non-physical oscillation of the numerical solution, unless $h \approx \varepsilon$ or  $h < \varepsilon$. 
Even if the mesh is adapted in the boundary layers regions, see e.g. \cite{CRD-results, Dahmen-Welper-Cohen12}, an SPLS discretization such as the  $P^1-P^2$ discretization can still produce non-physiscal oscillations. 

The goal of the paper is to analyze  the  presence of non-physical oscillations, 
 to propose ways to avoid such  oscillations, and to propose a new approach for establishing approximation errors  for convection dominated problems. 

The rest of the paper is organized as follows. 
Section \ref {sec:ON} contains a review of the main results on the optimal trial norms for the convection diffusion problem.  The Standard Linear (SL), the Saddle Point Least Square (SPLS) discretization  errors and non-physical oscillations analyses are presented in Section \ref{sec:Lin+SPLS}.  The one dimensional Upwinding Petrov-Galerkin (UPG) method emphasizing on non-oscillatory behavior of the discrete solution is summarized in Section \ref{sec:PG}.  In Section \ref{sec:UPGinfError}, for the one dimensional case, we motivate  the importance of recovering the exact solution at the nodes and the advantages of  discrete infinity norm approximation towards   eliminating non-physical oscillations.   In Section  \ref{ss:PG2D}, for a two dimensional case,  we review  the  quadratic bubble UPG  approximation properties and explain the oscillations along the parabolic boundary layer.  We summarize our findings  in Section \ref{sec:conclusion}.


 
 \section{Optimal trial norm for the convection diffusion problem}\label{sec:ON}
 We consider the variational formulation  \eqref{eq:var-cdr} 
 with $V=Q =H_0^1(\Omega)$   
and consider different  norms on the test and trial spaces.  On the test space  $V:=H_0^1(\Omega)$, we consider   the norm induced by 
 \[
a_0(u,v):=(\nabla u, \nabla v). 
\]
From the stability point of view, to define the optimal norm on $Q$,  we represent the {\it antisymmetric part}  of the bilinear form $b(\cdot, \cdot)$ in  the  $a_0(\cdot ,\cdot)$  inner product as follows.  Let $T:Q \to V$ be the representation operator defined by 
\[
a_0(Tu, v)  = (b \cdot \nabla u, v),  \Forall v\in V.
\]
It is easy to check  that
\[
|Tu| = \|b \cdot \nabla u\|_{H^{-1}(\Omega)} \leq \|b\| \|u\|_{L^2(\Omega)}.
\]

 For the one dimensional case and $\bb=1$, we have 
 \[
 a_0(Tu,q) = (u',q), \ \text{for all} \  q \in Q.
 \]
By solving the corresponding differential equation, one can find that
\begin{equation}\label{eq:T-Action}
Tu = x\overline{u} - \int_0^xu(s)\,ds, 
\end{equation}
and
\begin{equation}\label{eq:T-NormSq}
|Tu|^2 = \int_0^1 |u(s)-\overline{u}|^2\, ds= \|u -\overline{u}\|^2= \|u\|^2 -\overline{u}^2\leq \|u\|^2.
\end{equation}
Throughout the paper, for  $v\in L^2(0, 1)$, we denote   $\overline{v}:=\int_0^1 v(x)\, dx $.

 The continuous {\it optimal  trial norm} on $Q$ is defined by 
\[
 \|u\|_{*}:= \du {\sup} {v \in V}{} \ \frac {b(v, u)}{|v|} =   \du {\sup} {v \in V}{} \ \frac {\varepsilon a_0(u, v) + a_0(Tu,v) }{|v|}.
\]
Using the Riesz representation theorem and  $ a_0(Tu,u) =0$, we obtain  
\begin{equation}\label{eq:COptNormd}
\|u\|_{*}^2  =\varepsilon^2|u|^2 +|Tu|^2.
\end{equation}
For the one dimensional case,  from  \eqref{eq:T-NormSq}, we get 
\begin{equation}\label{eq:COptNorm}
\|u\|_{*}^2 = \varepsilon^2|u|^2 + \|u\|^2 -\overline{u}^2.
\end{equation}
Next, we assume  that $V_h\subset V=H_0^1(\Omega)$ and $\M_h\subset Q=H_0^1(\Omega)$  are discrete finite element spaces and that $\M_h\subset V_h$. 
To  describe  the {\it discrete optimal norm} on $\M_h$, we  will need  the standard elliptic projection  $P_h:Q\to V_h$ defined by 
\[
a_0(P_h\, u, v_h) = a_0(u,v_h), \ \text{for all} \, v_h \in V_h. 
\]
The optimal trial norm   on $ \M_h$ is 
 \begin{equation}\label{eq:dotn}
  \|u_h\|_{*,h}:= \du {\sup} {v_h \in V_h}{} \ \frac {b(v_h, u_h)}{|v_h|}.
 \end{equation}
As in the continuous case, see  \cite{connections4CD, comparison4CD},  by denoting $|u|_{*,h}:= |P_hTu|$, we have

\begin{equation}\label{eq:COptNorm-h}
\|u_h\|_{*,h}^2  =\varepsilon^2|u_h|^2 +|P_hTu_h|^2 = \varepsilon^2|u_h|^2 +|u_h|^2_{*,h}.
\end{equation}

The advantage of using  the  optimal trial norm on $Q$ and $\M_h$  resides with  the fact that both $inf-sup$ and $sup-sup$ constants  at the continuous and discrete levels are equal to one. 
The following error estimate  was proved in \cite{comparison4CD}. 
 \begin{theorem}\label{th:ap-PG2}
Let  $\|\cdot\|_{*}$ and $\|\cdot\|_{*,h}$ be the norms on $Q$, and $\M_h$ and  assume:
\begin{equation}\label{eq:c0}
\|v\|_* \leq c_0\|v\|_{*,h}\quad\quad\text{for all $v\in Q$}.
\end{equation}
 Let $u$ be the solution of \eqref{eq:var-cdr} and let  $u_h$ be  the unique solution of the problems \eqref{eq:var-cdr-h} or \eqref{SPLS4model-h}. Then, the following error estimate holds:
\begin{equation}\label{eq:Approx-Opt}
\|u-u_h\|_{*,h}\leq c_0\, \du{\inf}{p_h \in \M_h}{} \ \|u-p_h\|_{*,h}.
\end{equation}
\end{theorem}
 \subsection{Discrete optimal  trial norm for the one dimensional case} \label{sec:optimalnorm-CDc}
We review a representation  formula for $|\cdot|_{*,h}$  that was discussed in  \cite{CRD-results, connections4CD}. 

For  $V=Q=H^1_0(0,1)$ we consider the  standard inner product  given by  $a_0(u,v) = (u,v)_V = (u',v')$.
We  divide the interval $[0,1]$ into $n$  equal length subintervals using the nodes $0=x_0<x_1<\cdots < x_n=1$ and denote  $h:=x_j - x_{j-1}, j=1, 2, \cdots, n$. We define  the corresponding finite element discrete space  $\M_h$  as  the space of all {\it continuous piecewise linear  functions} with respect to the given nodes, that are zero at $x=0$ and $x=1$. 

Next, we let  $\M_h, V_h $  be the standard spaces of continuous piecewise linear functions 
\[
\M_h=V_h=span \{\varphi_1, \cdots,\varphi_{n-1}\}. 
\]
In this case,  the explicit formula for  $|u|_{*,h}= |P_hTu|$ is 
 \begin{equation}\label{eq:PhTu-Norm} 
|u|^2_{*,h}:= |P_hTu|^2 = \frac{1}{n}\sum_{i=1}^n \left(\frac{1}{h}\, \int_{x_{i-1}}^{x_i} u(x)\,dx\right)^2 - \left(\int_0^1u(x)\,dx\right)^2.
\end{equation}

\begin{remark}\label{rem:TS}
Note that $|\cdot |_{*,h}$ is a seminorm on $\M_h$. For $n=2m$ and  $\omega_h:=\varphi_1+\varphi_3+\cdots +\varphi_{2m-1}$,  we have $ |\omega_h|_{*,h}=0$. 
The graph of $\omega_h$  has a  {\it  ``teeth saw''  shape}  and can be highly oscillatory when $h=1/n$ is  small. 
\end{remark}

\section{The SL and SPLS discretizations exhibit non-physical oscillations}\label{sec:Lin+SPLS}
In this section, for the one dimensional case,  we  focus on  the Standard Linear (SL) discretization with $C^0-P^1$ test and trial spaces, and on  the SPLS discretization with $C^0-P^1$ trial space and $C^0-P^2$ test space. We  explain the oscillatory behavior for both discretizations  based on the error analysis in the optimal trial norms, and  on the  {\it closeness} between the discrete solution and the reated transport problems. 

For the finite element discretization, we use the following  notation:
\[ 
\begin{aligned}
a_0(u, v) & = \int_0^1 u'(x) v'(x) \, dx, \ (f, v) = \int_0^1  f(x) v(x) \, dx,\ \text{and}\\
b(v, u)& =\varepsilon\, a_0(u, v)+(u',v)  \ \text{for all} \ u,v \in V:=H^{1}_0(0,1).
\end{aligned}
 \]
A variational formulation of \eqref{eq:1d-model} is: Find $u \in V:= H_0^1(0,1)$ such that
 \begin{equation}\label{VF1d}
b(v,u) = (f, v), \ \text{for all} \ v \in V=H^{1}_0(0,1).
\end{equation}

\subsection{Standard discretization with $C^0-P^1$ test and trial spaces}\label{sec:1d-lin-discrete}

We  divide the interval $[0,1]$ into $n$  equal length subintervals using the nodes $0=x_0<x_1<\cdots < x_n=1$ and denote  $h:=x_j - x_{j-1}, j=1, 2, \cdots, n$. For the above uniform distributed notes on $[0, 1]$, we define  the corresponding finite element discrete space  $\M_h$  as  the subspace of $Q=H^1_0(0,1)$, given by
 \[ 
 \M_h = \{ v_h \in Q \mid v_h \text{ is linear on each } [x_j, x_{j + 1}]\},
 \]
  i.e., $\M_h$ is the space of all {\it continuous piecewise linear  functions} with respect to the given nodes, that are zero at $x=0$ and $x=1$.  We consider the nodal basis $\{ \varphi_j\}_{j = 1}^{n-1}$ with the standard defining property $\varphi_i(x_j ) = \delta_{ij}$.  
We couple the above discrete trial space with the  discrete  test space $V_h=\M_h$.  
 Thus, the standard $C^0-P^1$  variational formulation of \eqref{VF1d} is: \\ Find $u_h \in \M_h$ such that
 \begin{equation}\label{dVF}
b(v_h, u_h) =\varepsilon (u_h', v_h') + (u_h', w_h) = (f, v_h), \ \text{for all} \ v_h \in V_h.
\end{equation}
 In this case, according to Section \ref{sec:ON}, we have  
$
\|u_h\|_{*,h}^2  =\varepsilon^2|u_h|^2  +|u_h|^2_{*,h},
$
 where $|\cdot |^2_{*,h}$ has the representation given in \eqref{eq:PhTu-Norm}. 
 
As a consequence of  Theorem \ref{th:ap-PG2}, we proved in  \cite{CRD-results} the following result.
\begin{theorem}\label{th:opt-lin}
If $u$ is the solution of \eqref{VF1d} and $u_h$ is the solution of the  linear discretization \eqref{dVF}, then 
\[
\|u-u_h\|_{*,h}  \leq c_0 \du{\inf}{v_h \in V_h}{} \ \|u-v_h\|_{*,h}, \ \text{where} 
\]
\[
c_0=c(h,\varepsilon)= \sqrt {1+ \left (\frac{h}{\pi \, \varepsilon}\right )^2 } \approx \frac{h}{\pi \, \varepsilon} \ \text{ if }\ \varepsilon \ll h.
\]
\end{theorem}

 
  Numerical  tests  performed for \cite{CRD-results} show that, as $\varepsilon \ll h$, the linear finite element solution of \eqref{dVF}  presents non-physical oscillations.

To understand the presence of such oscillations, we consider the  reduced continuous and discrete  corresponding problems. By letting   $\varepsilon \to 0$  in \eqref{VF1d}, we obtain the {\it reduced continuous problem}: \\
Find $u \in H_0^1(0,1)$ such that
\begin{equation} \label{VF1d-simplified-w}
(u',v) = (f, v), \ \text{for all} \ v \in V. 
\end{equation}

The problem \eqref{VF1d-simplified-w} has unique solution,  if and only if $\int_0^1  f(x)  \, dx=0$.    On the other hand, the  corresponding Left to Right  {\it (LR) transport} problem: \\  Find $w \in H^1(0,1)$ such that
 \begin{equation}\label{VF1d-reduced}
w'(x) = f(x)  \ \text{for all} \ x \in (0, 1), \text{and} \ w(0)=0,
\end{equation}
has  unique solution: $w(x) = \int_0^x  f(s)  \, ds$, regardless of the average of the function $f$. 

Similarly,  the other related Right to Left  {\it  (RL) transport} problem: Find $\theta \in H^1(0,1)$ such that
 \begin{equation}\label{VF1d-reduced2}
\theta'(x) = f(x)  \ \text{for all} \ x \in (0, 1), \text{and} \ \theta(1)=0, 
\end{equation}
has  unique solution: $\theta(x)= w(x) -\int_0^1 f(x)\, dx$.
 
 By letting   $\varepsilon \to 0$  in \eqref{dVF}, we obtain the {\it  reduced discrete problem}: \\
Find $U_h \in \M_h$ such that
\begin{equation} \label{VF1d-simplified-wh}
(U_h',v_h) = (f, v_h), \ \text{for all} \ v_h \in V_h=\M_h. 
\end{equation}

\hspace{-12mm} \parbox{0.55\textwidth}{
\begin{center}
\includegraphics[width=0.55\textwidth]{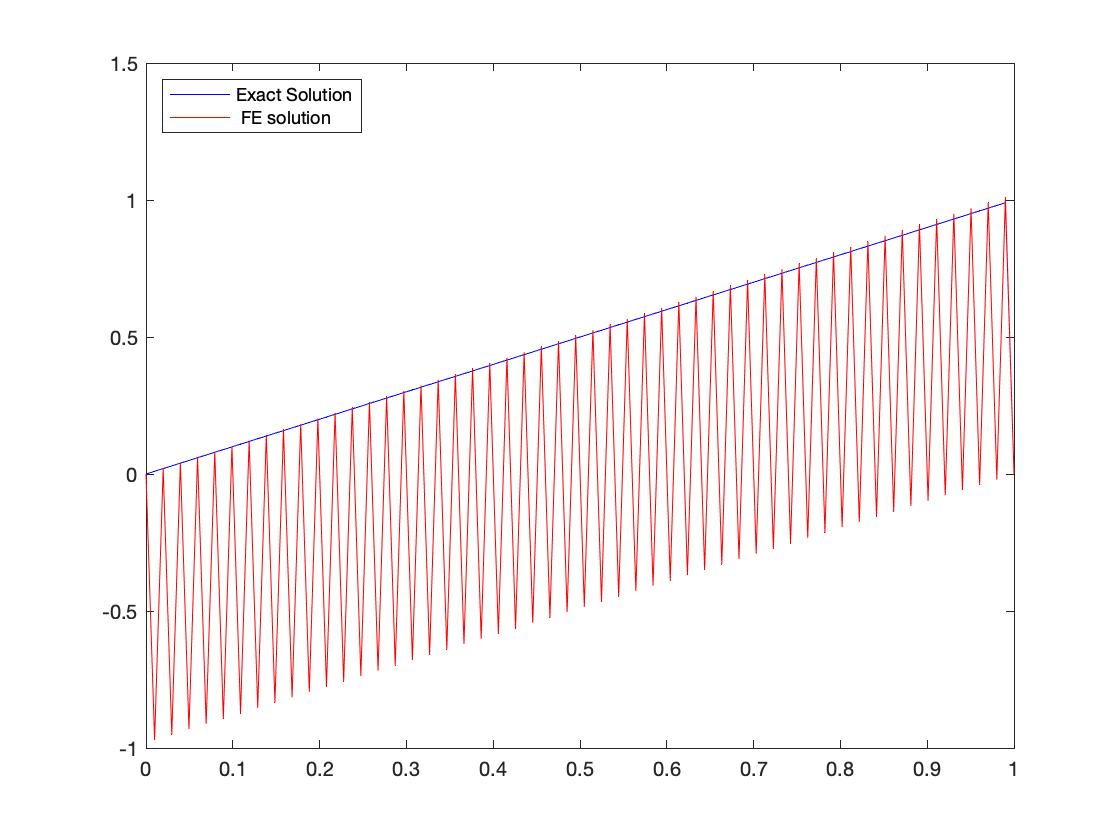}\\
{Fig.1: $f=1, n=99, \varepsilon=10^{-6}$ \\ $u_h$ oscillate {between $x$ and $x-1$} and \\ very close to the solution $U_h$ of \eqref{VF1d-simplified-wh}} 
\label{fig:Fig1}
\end{center}
}
\parbox{0.55\textwidth}{
\begin{center}
\includegraphics[width=0.55\textwidth]{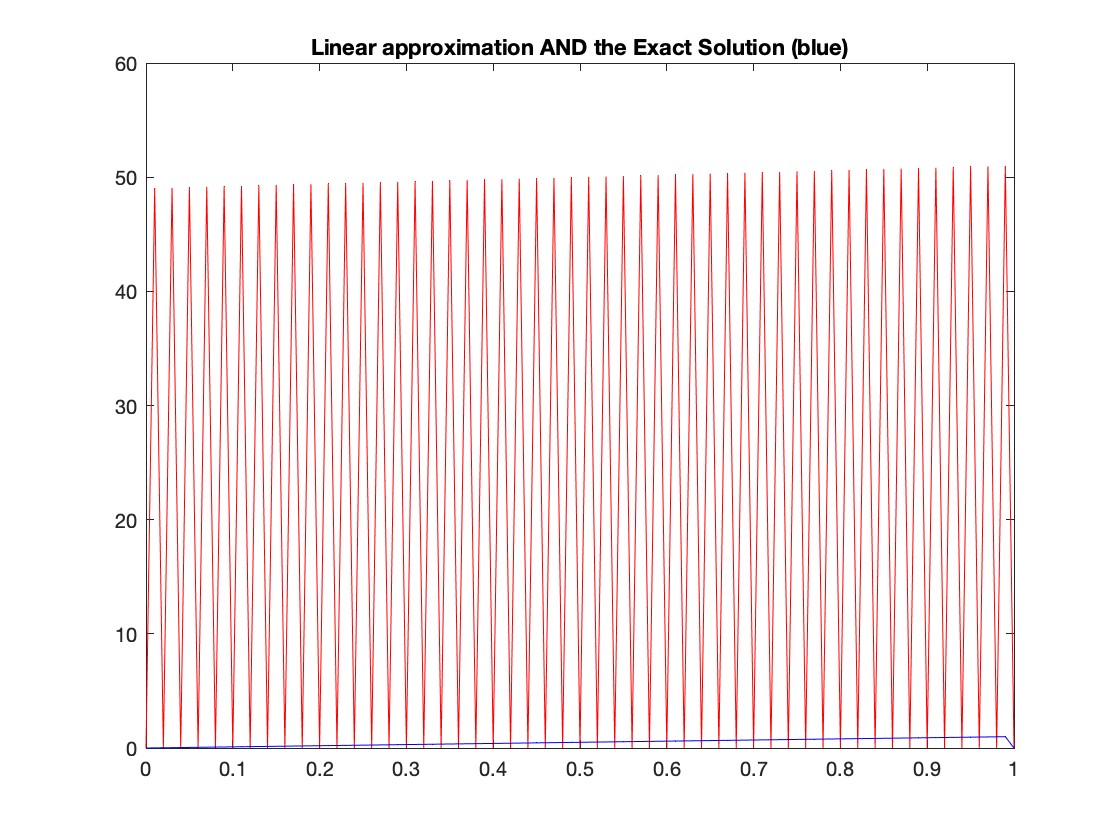}\\
{Fig.2: $f=1, n=100$,  $\varepsilon=10^{-6}$ \\ $u_h  \approx  u + h^2/(2\varepsilon) \  \omega_h$. {\it No solution for the  reduced discrete  problem} \eqref{VF1d-simplified-wh}.} 
\label{fig:Fig2}
\end{center}
}

 If $\int_0^1 f(x) \, dx\neq 0$ and the number of subintervals $n=1/h$  is odd, then the linear system associated with \eqref{VF1d-simplified-wh} has  unique solution, see  \cite{CRD-results}. In this case,  as $\varepsilon/h \to 0$, the finite element solution $u_h$ of \eqref{dVF} is very close to  the solution  $U_h$ of  \eqref{VF1d-simplified-wh}, and both solutions  oscillate between the graph of the two transport problems RL and LR, i.e.,   $w$  and $\theta$, respectively.
 
 If $\int_0^1  f(x)  \, dx\neq 0$ and the number of subintervals $n$  is even,  the problem  \eqref{VF1d-simplified-wh}  might not have solution, but because of the week coercivity of the bilinear form $b(\cdot, \cdot) $ of \eqref{dVF}, we have  that \eqref{dVF} has an unique solution $u_h$. 
 Numerically, as shown in Figure 2, we observed  that for  ${\varepsilon}/{h} \leq 10^{-4}$, at the nodes,  $u_h$ is very close to  $u + h^2/(2\varepsilon)\,  \omega_h$, where $u$ is the  solution of \eqref{VF1d}  and $\omega_h$ is the  “teeth saw'' function defined in Remark \ref{rem:TS}. Thus, at the nodes, we have that $u_h-u$  behaves like $h^2/(2\varepsilon)\ \omega_h$.  

Regardless of the number of intervals,  for $\int_0^1  f(x)  \, dx\neq 0$ and  for $\varepsilon<h$, for example  
$\varepsilon/h <1/3$,  we observe non-physical oscillations in the SL discrete solution. 

If the  {\it  reduced discrete problem} does not have a solution, we can expect global non-physical oscillatory behavior of the discrete solution $u_h$  of  \eqref{dVF}, and $u_h$  could contain high frequency modes such as the  “teeth saw'' function.

 If the {\it  reduced discrete problem} does  have a unique solution, but it does not correspond to the discretization of a reduced  {\it  continuous problem},  we can also  expect global non-physical  oscillations  in the discrete solution. 

Theorem \ref{th:opt-lin}  cannot control   the discrete infinity error of the solution, hence it cannot predict the  oscillatory behavior of the discrete solution.




\subsection{Approximation for the SPLS discretization} \label{sec:SPLS}
A  {\it  saddle point least square}  (SPLS) approach   for solving \eqref{eq:1d-model} has  been used before, for example in \cite{CRD-results, BM12, Dahmen-Welper-Cohen12}.  For $V=Q= H^1_0(0,1)$, we look for  finding $(w=0, u) \in V \times Q$ such that 
\begin{equation}\label{SPLS4model2}
\begin{array}{lclll}
a_0(w,v) & + & b( v, u) &= (f,v ) &\ \Forall  v \in V,\\
b(w,q) & & & =0   &\  \Forall  q \in Q,   
\end{array}
\end{equation}
where
\[
b(v, u) =\varepsilon\, a_0(u, v)+(u',v)= \varepsilon\, (u', v')+(u',v). 
\]
Clearly, the  component $u$ of the solution $(w=0, u)$ is the solution of  \eqref{VF1d}. The advantage of considering the SPLS form is that its discretization leads to a symmetric  linear system. 
Analyses for finite element  test and trail spaces of various degree polynomials for \eqref{SPLS4model2} can be found in \cite{Dahmen-Welper-Cohen12, dem-fuh-heu-tia19}. 

For $P^1-P^2$ discretization we follow the notation of \cite{CRD-results}. We consider on $[0, 1]$ the uniformly distributed nodes $x_j=hj, j=0,1,\cdots, n$, with  $h=1/n$ and  define $\M_h= C^0-P^1:= span\{ \varphi_j\}_{j = 1}^{n-1}$, with $\varphi_j$'s  the standard linear nodal functions  and $V_h=C^0-P^2$- the space of continuous piece-wise quadratic functions.
The discrete version of \eqref{SPLS4model2} is: \\ 
 Find $(w_h, u_h) \in V_h \times \M_h$ such that \eqref{SPLS4model-h}
is satisfied.  
The component $u_h$ of  $(w_h, u_h)$ is the SPLS discretization of  \eqref{eq:1d-model} or the SPLS solution of  \eqref{SPLS4model-h}.

We note that the projection $P_h$  defined in Section \ref{sec:ON}, is the projection on the space $V_h=C^0-P^2$ of continuous piece-wise quadratic functions. For any piecewise linear function $u_h \in \M_h$, we have that 
\[
Tu_h = x\overline{u}_h - \int_0^xu_h(s)\,ds,
\]
is a continuous piecewise quadratic function which is zero at the ends. Consequently, $Tu_h \in V_h$, and $P_h\, Tu_h =T u_h$,  
and, according to  \eqref{eq:COptNorm-h} and \eqref{eq:COptNormd}, the optimal discrete norm on $\M_h$ becomes 
\[
\|u_h\|_{*,h}^2 =\varepsilon^2 |u_h|^2 + |Tu_h|^2 =\|u_h\|_{*}^2.
\]
Using the optimal norm on $\M_h$, a discrete $\inf-\sup$ condition is satisfied, and the problem \eqref{SPLS4model-h} has a unique solution. 
In addition, for  the  $P^1-P^2$ SPLS discretization, we can consider the same norm given by 
\[
\|u\|_{*}^2 =\varepsilon^2|u|^2 +\|u-\overline{u}\|^2=\varepsilon^2|u|^2 +\|u\|^2-\overline{u}^2=\|u\|_{*,h}^2 
\]
 on both spaces $Q$ and $\M_h$. Since \eqref{eq:c0} is satisfied with $c_0=1$,  as a consequence of Theorem \ref{th:ap-PG2},  we obtain:
\begin{theorem}\label{th:P1P2err}
If $u$ is the solution of \eqref{VF1d}, and $u_h$ is the  SPLS solution  of  the $(P^1-P^2)$ discretization \eqref{SPLS4model-h}, then  
$$
\|u - u_h\|_{*}\leq \inf_{p_h\in \M_h}\|u - p_h\|_{*} \leq \|u - u_I\|_{*},
$$
where $u_I$ is the interpolant of the exact solution on the $h$-uniformly distributed nodes on $[0, 1]$. 
\end{theorem}

For  $\int_0^1f(x)\, dx=0$, the $P^1-P^2$ SPLS discretization improves on standard linear discretization of \eqref{VF1d} from the error point of view, see e.g. \cite{CRD-results, Dahmen-Welper-Cohen12}. For  $\int_0^1  f(x) \, dx\neq 0$ and  $\varepsilon \ll h$, in spite of the optimal  approximation result of Theorem \ref{th:P1P2err}, the SPLS solution $u_h$ approximates a shift by a constant of the solution $u$ of the continuous problem \eqref{SPLS4model2} or \eqref{VF1d}, and  {\it  local non-physical oscillations} still appear in the discrete solution  $u_h$ at the ends of the interval. More numerical experiments showing the oscillations are in \cite{CRD-results}. 

\subsection{The oscillatory behavior of  $P^1-P^2$ SPLS discretization} \label{sec:behaveSPLS}
 

 In this section, we justify why the $P^1-P^2$ SPLS discretization fails to produce a good approximation for the solution of of \eqref{VF1d} for the case $\int_0^1  f(x) \, dx\neq 0$. In what follows, we use the notation and considerations of Section \ref{sec:SPLS}. 
 
We consider  the {\it  reduced continuous problem} obtained from \eqref{SPLS4model2} by letting $\varepsilon \to 0$, i.e.,  Find $(w, u) \in V \times Q$ such that 
\begin{equation}\label{SPLS4model-h-R}
\begin{array}{lclll}
(w',v') & + & (u', v) &= (f,v ) &\ \Forall  v \in V= H^1_0(0,1)\\
(w, q') & & & =0   &\  \Forall  q \in Q = H^1_0(0,1).  
\end{array}
\end{equation}
The problem  is not well posed when $\int_0^1  f(x)  \, dx\neq 0$. We can change the trial space $Q$ to 
$L^2_0(0, 1):=\{u \in L^2(0, 1) | \ \int_0^1 u=0\}$ in order to have  existence and uniqueness of the solution. Nevertheless,  in this case, the solution space cannot satisfy  the boundary conditions of the original problem \eqref{eq:1d-model}.  

However, the  {\it reduced discrete problem} obtained from \eqref{SPLS4model-h}  by letting \\ $\varepsilon \to 0$, i.e., Find $(w_h, u_h) \in V_h \times \M_h$ such that 
\begin{equation}\label{SPLS4model-h-R}
\begin{array}{lclll}
(w'_h,v'_h) & + & (u'_h, v_h) &= (f,v_h ) &\ \Forall  v_h \in V_h,\\
(w_h,q'_h) & & & =0   &\  \Forall  q_h \in \M_h,   
\end{array}
\end{equation}
 has  unique solution. This is justified by the fact that a discrete $\inf-\sup$ condition, using optimal trial norm  on $\M_h$ and the standard $H^1$ seminorm on $V_h$, holds.  
 
 Numerical tests in \cite{CRD-results} showed that oscillation of  {the discrete solution}  $u_h$ of \eqref{SPLS4model-h-R} predict oscillatory behavior of  the SPLS discrete solution $u_h$ of \eqref{SPLS4model-h}. In fact, for $\varepsilon/h \leq 10^{-4}$, the two solutions look identical in the ``eye ball measure".
 
 Next,  we explain  the solution   behavior for  the {\it reduced discrete problem} \eqref{SPLS4model-h-R}. 
We introduce   $u^f \in V=H_0^1(0,1)$ as the solution of
 \begin{equation}\label{eq:uf} 
 -(u^f)^{''}= f, \ 0<x<1 \ \ \text{or}   \ ((u^f)',v') = (f,v ),  \ \Forall  v \in V, 
 \end{equation}
and the elliptic projection of $u^f$ on $V_h= C^0-P^2$, as  the solution $u_h^f \in V_h$ of
 \begin{equation}\label{eq:ufh} 
 ((u^f_h)',v'_h) = (f,v_h ),  \ \Forall  v_h \in V_h.
 \end{equation}
 In this section,  we  also need the solution $w(x) = \int_0^x  f(s)  \, ds$  of the  {\it (LR) transport} problem \eqref{VF1d-reduced}, and the following two subspaces of $L^2(0, 1)$:
 \[
\overline{\M}_h:= \{w_h -\overline{w}_h \ | w_h \in \M_h\}, \ \text{and}
\] 
 \[
\tilde{\M}_h:=\overline{\M}_h \bigoplus  span \{1\} =\{v_h\in C^0-P^1 \ | \ v_h(0)=v_h(1)\}.
\]
Next, we  state the main result of this section.
\begin{theorem} \label{th:SPLSosc}
Let $u_h$ be the solution of  the {\it reduced discrete problem} \eqref{SPLS4model-h-R}. Then,  
{ $u_h-\overline{u}_h$ is the $L^2(0,1)$ orthogonal projection of $w(x) -\overline{w} $ onto $\tilde{\M}_h$}.
\end{theorem}
\begin{proof}

 The system \eqref{SPLS4model-h-R} is equivalent to 
\begin{equation}\label{SPLS4model-h-R2}
\begin{aligned}
w_h  +  P_h T u_h &= u^f_h , \\
(w_h',q_h)  \ \ \ \ \ & =0   \  \Forall  q_h \in \M_h,
\end{aligned}
\end{equation}
where $\M_h =span \{\varphi_1, \cdots,   \varphi_{n-1}\}$, and $P_h$ is the elliptic ptojection on $V_h$.

As presented in Section \ref{sec:SPLS},  $P_h T u_h=T u_h$ and $(T u_h)' = \overline{u}_h - u_h$. 
Differentiating  the first equation of \eqref{SPLS4model-h-R2}, we obtain 
\begin{equation}\label{SPLS4model-h-R3}
w_h'  + \overline{u}_h - u_h = (u^f_h)'.  
\end{equation}
From the second  equation of \eqref{SPLS4model-h-R2},  using that $\int_0^1 w_h' =0$, we have
\begin{equation}\label{eq:whprime}
 (w_h',q_h - \overline{q}_h) =0   \  \Forall  q_h \in \M_h.
 \end{equation}

By substituting  $w_h'$ from \eqref{SPLS4model-h-R3}  in \eqref{eq:whprime}, we obtain  
\begin{equation}\label{SPLS4model-h-R4}
 (u_h-\overline{u}_h, q_h - \overline{q}_h) =(-(u^f_h)' , q_h - \overline{q}_h),  \  \Forall  q_h \in \M_h.
\end{equation}
 Hence, $u_h-\overline{u}_h$ is the $L^2(0,1)$ orthogonal projection of $-(u^f_h)'$ onto $\overline{\M}_h$. 
 
 Next, we  prove that $u_h-\overline{u}_h$ is the $L^2(0,1)$ orthogonal projection of $-(u^f)'$ onto $\overline{\M}_h$. We note that for any $q_h - \overline{q}_h \in \overline{\M}_h$,  the function 
\[
v_h = \int_0^x (q_h(s) - \overline{q}_h)\, ds \in V_h= C^0-P^2, 
\]
 and 
  \begin{equation}\label{SPLS4model-h-R5}
 \begin{aligned}
((u^f_h)' , q_h - \overline{q}_h) & =  ((u^f_h)' , (v_h)') =(f, v_h) \\ &= ((u^f)' , (v_h)')=((u^f)' , q_h - \overline{q}_h).
\end{aligned}
\end{equation}
 From \eqref{SPLS4model-h-R4} and \eqref{SPLS4model-h-R5}, we have 
 \[
  (u_h-\overline{u}_h, q_h - \overline{q}_h) =(-(u^f)' , q_h - \overline{q}_h), \ \text{for all} \  q_h \in \M_h.
 \]
Thus,  $u_h-\overline{u}_h$ is indeed  the $L^2$ orthogonal projection of $-(u^f)'$ onto $\overline{\M}_h$.  
In addition, since $(u^f)'$ is $L^2$ orthogonal to  constant functions,  we also have that 
{ $u_h-\overline{u}_h$ is the $L^2$ orthogonal projection of $-(u^f)'$ onto $\tilde{\M}_h$}.

 Integrating $-(u^f)''=f$  on $[0, x]$ gives us 
  \begin{equation*}\label{eq:ufprime}
  { -(u^f)'(x) = \int_0^x f(s)\, ds -(u^f)'(0) :=w(x) -(u^f)'(0)}.
\end{equation*}
 Integrating  the above identity on $[0, 1]$ and using $(u^f)(0) = (u^f)(1) =0$  give
  \begin{equation*}\label{eq:ufprime0}
(u^f)'(0) =  \int_0^1 w(x)\, dx. 
\end{equation*}
From  the last two identities, we obtain 
\[
-(u^f)'(x) =w(x) - \overline{w},
\] 
 and  the theorem is proved. 
 \end{proof}
 Now, we  justify the behavior of the solution $u_h$ of the {\it reduced discrete problem} \eqref{SPLS4model-h-R}: 
 When  { $\int_0^1 f(x)\, dx \neq 0$},  the function { $u_h-\overline{u}_h$ is the   $L^2$ projection of the  continuous function $w(x) - \overline{w}$   {\it with different end values} to a space  $\tilde{\M}_h $ of continuous  functions on $[0, 1]$ with {\it the same values at the ends}. 
 For $f=1$, we get $u^f=0.5 (x-x^2)$ and $w(x) - \overline{w} =x-1/2$ which takes the values $\pm \frac12$ at $x=0$ and $x=1$, respectively. Since functions in  $\tilde{\M}_h$ have the same values at the ends, the $L^2$ orthogonal projection of $x-1/2$  onto $\tilde{\M}_h$ oscillates at the ends, see Figure \ref{fig:FigSPLsosc} (Left).
 
 \begin{figure}[h]
	\begin{center}
	\includegraphics[width=1.6in]{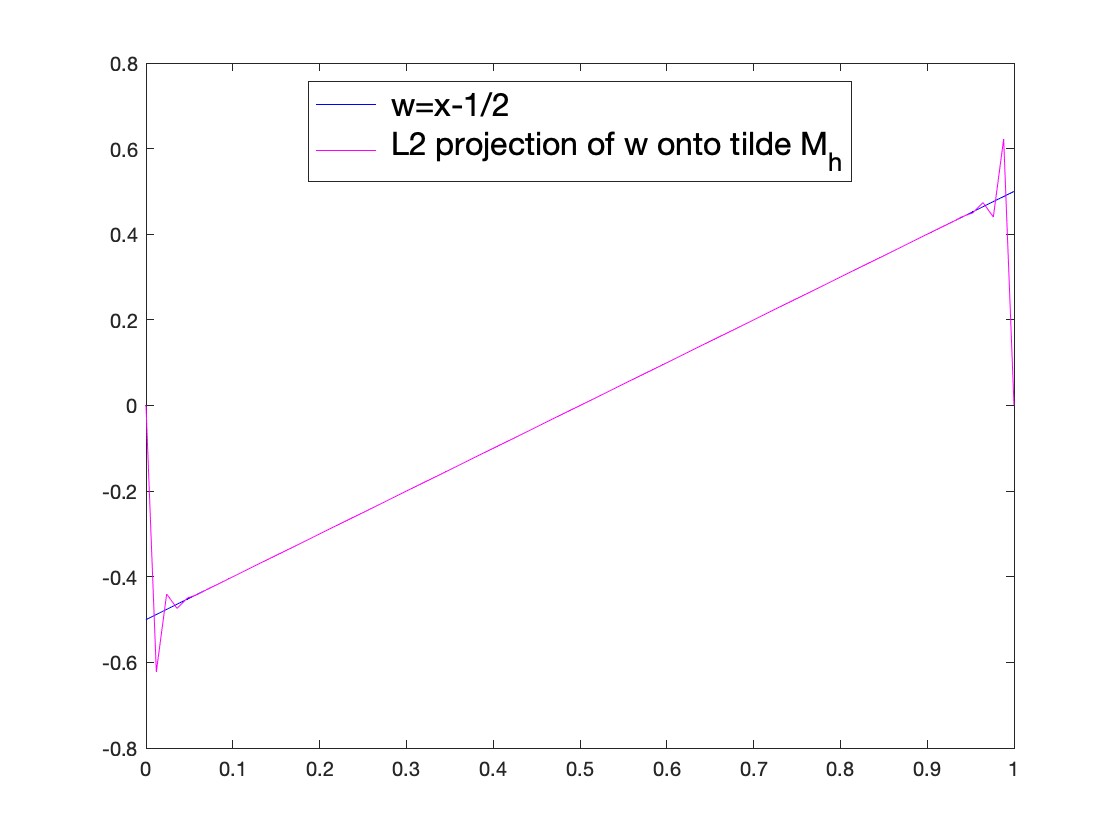}
\includegraphics[width=1.6in]{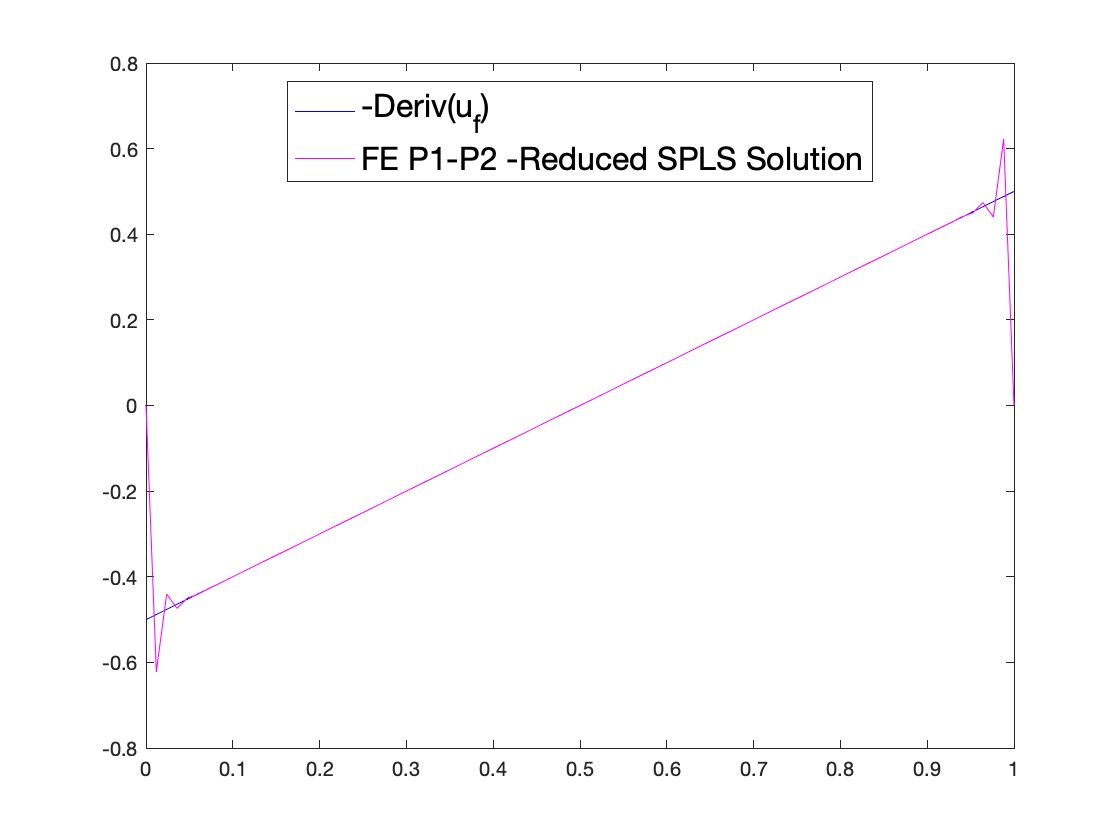}
\includegraphics[width=1.6in]{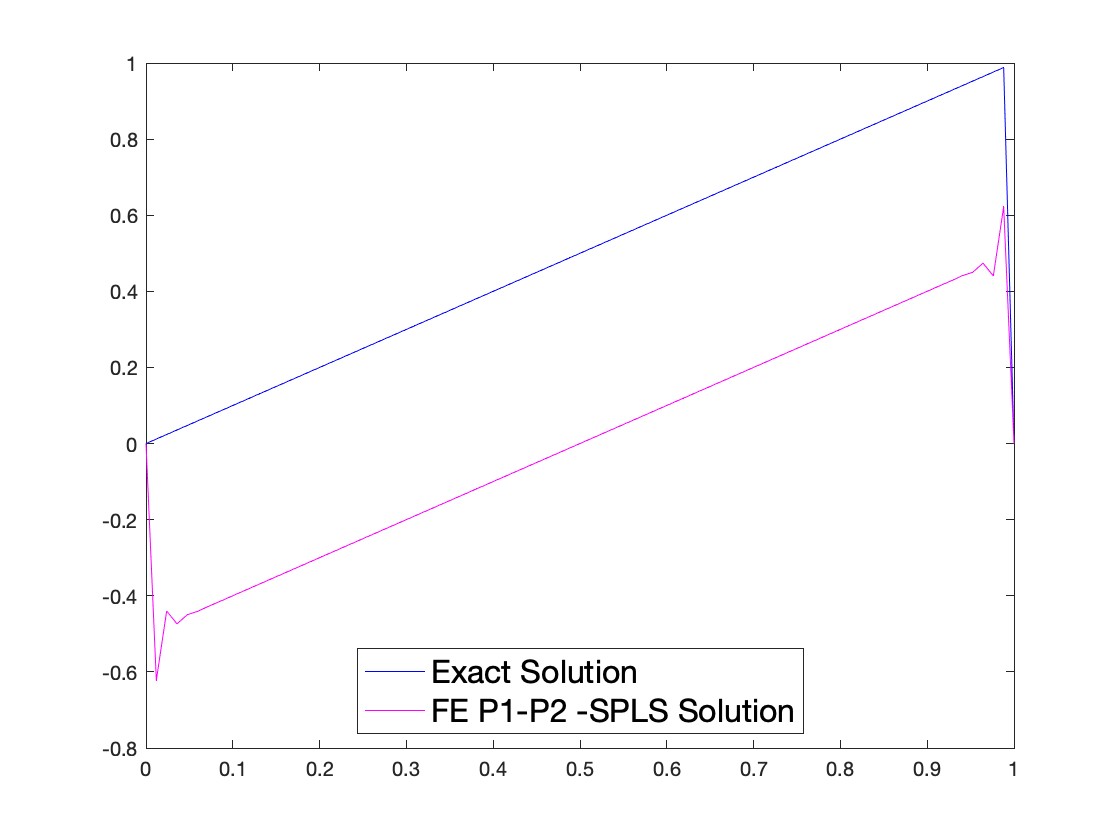}
\caption{$n=84$,  Left: $L^2$ proj. of $w=x-1/2$ on $\tilde{\M}_h$ \\ 
Middle: $P^1-P^2$  SPLS for Reduced Problem \\
Right:  $P^1-P^2$  SPLS for the Standard Problem, $\varepsilon=10^{-6}$} 
\label{fig:FigSPLsosc}
\end{center}	
\end{figure}
 
  When comparing  the solution of the SL  discretization to the  discrete solution $u_h$ of the standard SPLS formulation, we note that the discrete solution  $u_h$ is free of {\it global non-physical oscillations},  but exhibits  {\it local oscillations}  at the ends of the interval $[0, 1]$,  see Figure \ref{fig:Fig1} (Right). Numerical tests   for the case $\int_0^1  f(x) \, dx\neq 0$ and  $\varepsilon \ll h$ show that the SPLS discrete solution is very close to the solution of the {\it reduced discrete problem} \eqref{SPLS4model-h-R}. Both solutions, except for the end oscillations, are  very close to the graph of the function 
   \[
 {U(x):=w(x) -\frac12\, w(1)= \int_0^x f- \frac12 \int_0^1 f }={1/2(w(x)  +\theta(x))}, 
\]
where $w$ and $\theta$ are the solutions of the  the RL and LR transport problems.


\begin{remark} \label{rem:SPLS1d}
 {\it The reduced problem} \eqref{SPLS4model-h-R} is independent of $\varepsilon$, and its solution can still exhibit oscillations.    For  $\varepsilon \ll h$, the reduced problem  \eqref{SPLS4model-h-R} predicts oscillatory behavior of the  solution of the standard SPLS discretization \eqref{SPLS4model-h}.  Oscillations at the boundary of the domain,   appear  because the discrete solution is approximated  by an $L^2$- projection of a continuous  function with non-zero boundary conditions to a  subspace of continuous functions  that does not  account for the  boundary conditions of the function. 
  \end{remark}
  
Both  SL and SPLS reduced discrete problems  do not correspond to continuous reduced problems  that have unique solutions. 

For both SL and SPLS discretizations, the non-physical oscillations are related with the LR and RL transport problems. 
It has been proven in \cite{connections4CD} that for $\varepsilon \ll h$, the exact solution $u$ of \eqref {eq:1d-model}  and the solution $w$  of the LR transport  problem are very close at the interior nodes $x_1, x_2, \cdots, x_{n-1}$. We will address this idea again in Section \ref{asec:ConvDiscInf}. Thus, for any finite element discretization of \eqref {eq:1d-model}, the closeness of the discrete solution $u_h$ to  $w$ can lead to an oscillation free discretization. However, the ``closeness'' of $u_h$ to the solution $\theta$ of the RL transport  problem can lead to non-physical oscillations. 

 
\section{The Upwinding Petrov-Galerkin method  with bubble type test space } \label{sec:PG}
In this section, we review the one dimensional UPG method emphasizing on non-oscillatory behavior of the discrete solution. Our  main idea of the bubble  UPG method  is to  {\it  choose the test space} to {\it create upwinding diffusion from the convection part}. To create a basis for the test space, for each nodal basis of the trial space, we add  locally supported upwinding bubbles. The {\it upwinding process}  leads to  the elimination of the non-physical oscillation in the discrete solutions, and to better approximation.  

To be more precise, we define the test space $V_h$ by introducing a  bubble function for each interval $[x_{i-1}, x_i], i=1,2, \cdots, n$. First, consider a {\it bubble generating function} $B:[0,h] \to\R$ with the following properties:
\[
B(0)=B(h)=0,\ \  \frac{1}{h} \int_0^h B ={b}, \ \text{and}  \int_0^h (B'(x))^2  \, dx=\frac{b_e}{h}\ \text{with} \ {b, b_e>0}. 
\]
Next,  for $ i=1, \cdots, n$,  define $B_i:[0, 1] \to \R$ by  {$B_i(x)=B(x-x_{i-1})$} on $[x_{i-1}, x_i]$,  and  extend it  by zero to the entire interval $[0, 1]$. The discrete test space  for the  bubble UPG discretization is 
\[
V_h:=  span \{ g_i: \, i=1,2,\cdots,n-1\}, \ \text{where} \ g_i= \varphi_i + ( B_{i}-B_{i+1}).
\]
The idea of  bubble enriched test spaces  was introduced  for quadratic bubbles in \cite{zienkiewicz76}. A generalization of the method, using an arbitrary  generating bubble function, is presented and analyzed in \cite{CRD-arxiv, CRD-results}. 
A particular choice for the generating bubble $B$  is  {$B=4\,  \varphi_{0} \varphi_{1}= 4 \,  \frac{x}{h}\left (1-\frac{x}{h}\right)$}. In this case, the graphs of a basis trial function $\varphi_i$ and   a basis test function $g_i=\varphi_i + (B_i-B_{i+1})$ are  presented below.  

\begin{figure}[h!]
\begin{center}
	\includegraphics[width=0.3\textwidth]{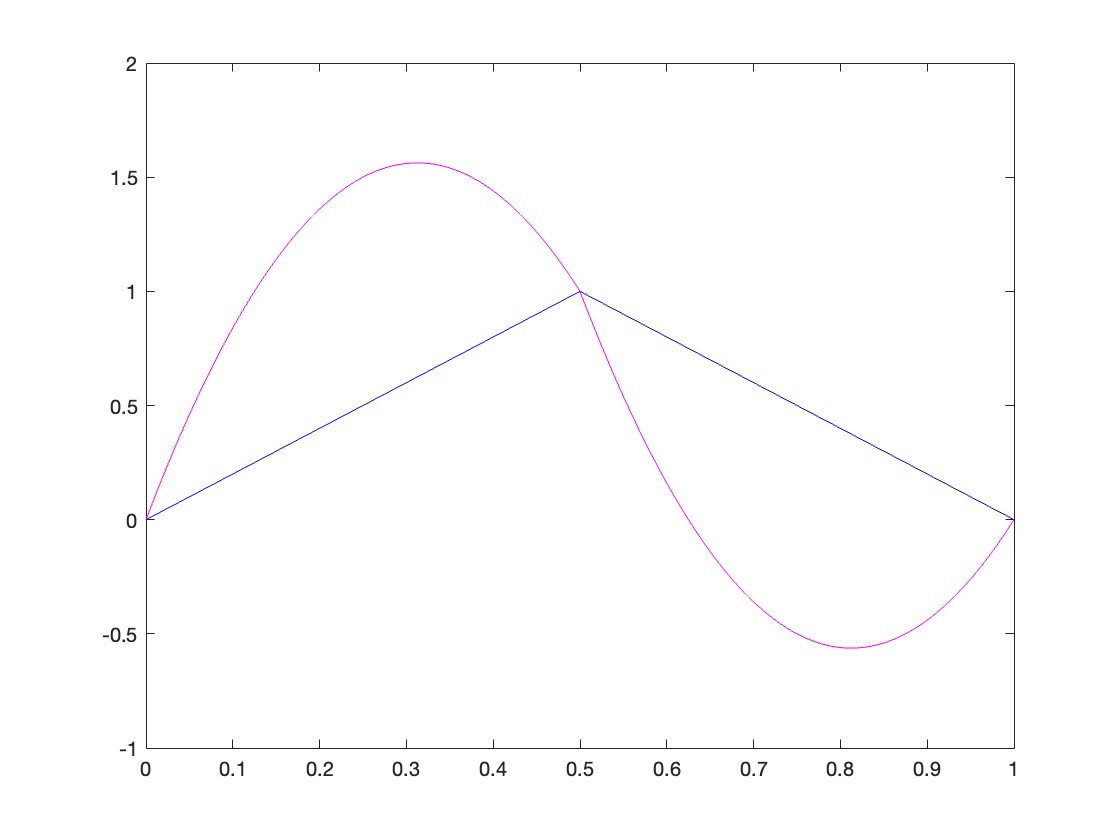}
\end{center}
\end{figure}
{The bubble UPG idea} for  solving \eqref{eq:1d-model} 
resides in the choice of the signs of the added bubbles in order to match the convection coefficient ${1>0}$. 

The upwinding Petrov Galerkin discretization with general bubble functions for 
\eqref{eq:1d-model} is: Find $u_h \in \M_h$ such that 
\begin{equation}\label{eq:1d-modelPG}
b(v_h, u_h) = \varepsilon\, a_0(u_h, v_h)+(u'_h,v_h) =(f,v_h) \ \Forall  v_h \in V_h. 
\end{equation}
As presented in \cite{connections4CD, comparison4CD}, by  writing {\it a generic test function  $v_h \in V_h$ as $v_h=w_h + B_h \in V_h$}, with {$w_h \in \M_h$} the linear part of $v_h$, and $B_h$ the bubble part of $v_h$, the problem \eqref{eq:1d-modelPG} with $V_h:=  span \{ g_i\}= span \{\varphi_i + ( B_{i}-B_{i+1})\}$ has the  reformulation:  Find $u_h \in \M_h=C^0-P^1$ such that 
\begin{equation}\label{eq:UPGvf}
\begin{aligned}
 \left (\varepsilon +b\, h \right )  (u'_h, w'_h) +  (u'_h, w_h) & =F_{upg}(w_h), \ \text{where}\\  
F_{upg}(w_h) &:=(f, w_h )  + h\,  (f, w'_h  \sum_{i=1}^{n}   B_i),  w_h \in \M_h. 
\end{aligned}
\end{equation}
Note that the  bubble part $B_h$ in the  convection term of the original variational formulation, i.e.,  $(u'_h, B_h)$,  produces  the extra diffusion  term  $b\,  h (u'_h, w'_h)$. In other words, the bubble part of the test space creates {\it upwinding diffusion} from the convection term. This  is  an important  feature of the bubble UPG method that can be extended to the multidimensional case.  

The reformulation \eqref {eq:UPGvf} involves only linear functions in the bilinear form.  
From \eqref {eq:UPGvf}, it is easy to check that the matrix associated with the reformulation is tridiagonal.  More precisely, 
\begin{equation} \label{eq:M4CDfe} 
M_{fe}= tridiag \left ( -{\left(\frac{\varepsilon}{h} +b\right)} - \frac{1}{2},\  2 {\left(\frac{\varepsilon}{h} +b\right)},
\   -{\left(\frac{\varepsilon}{h} +b\right)}+ \frac{1}{2} \right ), 
\end{equation}
and depends only on $\varepsilon, h$ and the average  $b$ of the generating bubble $B$. 

Based on the reformulation \eqref {eq:UPGvf}, we proved  the following error estimate  in the discrete  optimal trial norm, in  \cite{connections4CD, comparison4CD}. 

\begin{theorem}\label{thm:PGError}
If $u$ is the solution of \eqref{VF1d}, $u_h$ is the solution of the upwinding PG formulation \eqref{eq:1d-modelPG},  and the bubble  $B$ is chosen such that  $b\geq 1/\pi$, then
 
\begin{equation}\label{eq:PGerrB}
\|u - u_h\|_{*,h}\leq \sqrt{1+b_e}\, \inf_{p_h\in \M_h}\|u - p_h\|_{*,h} \leq \sqrt{1+b_e}\, \|u - I_h(u)\|_{*,h}, 
\end{equation}
where   $I_h(u)$ is the linear interpolant of $u$ on the uniform mesh of size $h$, 
\[
\|u_h\|_{*,h}^2  = \frac{(\varepsilon +h\, b)^2}{1+b_e}\, |u_h|^2 +  \frac{1}{1+b_e}\ |u_h|^2_{*,h}, 
\]
and   $|\cdot |_{*,h}$  is the seminorm defined in \eqref{eq:PhTu-Norm}. 
\end{theorem}
Note that  the ``teeth saw''  function $\omega_h$ satisfies  $|\omega_h|^2_{*,h}=0,\  |\omega_h|^2 =1/h^2$.  Hence,   
\[
\|\omega_h\|_{*,h}^2  =\frac{1}{1+b_e}  \frac{(\varepsilon +h\, b)^2}{h^2} \geq \frac{b^2}{1+b_e}\ \text{and}
\]
 the discrete solution  could not contain  spurious modes  of ``teeth saw'' type.  

The error estimate  \eqref{eq:PGerrB} is  quasi-optimal  in the discrete optimal norm, but it essentially depends 
on the approximation properties of the linear interpolant on uniform meshes. Since, in general, the solution $u$ exhibits a boundary layer close to $x=1$, the interpolant of $u$ fails to approximate well the exact solution  on regions close to $x=1$. Thus, the global  estimate \eqref{eq:PGerrB} might not  be able to control  the discrete infinity error, and therefore, the oscillatory behavior of the discrete solution. 

However,  for  special choices of the generating bubble function $B$, we can control the discrete infinity error of the bubble UPG method. For example, we can   choose a  bubble $B$  with the  average  $b=\frac12 -\frac{\varepsilon}{h}$, 
such that  the upper diagonal  entries in the matrix $M_{fe}$ are zero. The matrix $M=M_{fe} $  becomes: 
\begin{equation} \label{eq:M4CDut}
M= tridiag  ( -1,  1, 0), 
\end{equation}
and the linear system corresponding to \eqref{eq:UPGvf}  is 
 \[
 M\, U= F_{pg}, \ \text{where} \ F_{pg}=[(f,g_1), (f,g_2),\cdots, (f,g_{n-1})]^T.    
 \]
A  forward solve  leads  to the following formula for the component $u_j$ of $U$: 
\begin{equation}\label{eq:uj}
u_j= (f, \varphi_1 + \varphi_2+ \cdots + \varphi_j) + (f, B_1-B_{j+1}), \ j=1,2,\cdots, n-1.
\end{equation}
Introducing  the nodal function $\varphi_0$ corresponding to $x_0=0$, i.e., $\varphi_0$ is the continuous piecewise linear function on $[0, 1]$ such that $\varphi_0(x_j) =\delta_{0,j} $, \\ $ j=0,1,\cdots, n$, and  using that $\varphi_0+\varphi_1+\cdots +\varphi_j =1$ on $[0, x_j]$, the formula 
\eqref{eq:uj} leads to 

\begin{equation}\label{eq:eplpicitSol}
u_j= \int_0^{x_j} f(x)\ dx\,  + \int_0^{x_1} f(B_1-\varphi_{0})\ dx\,  + \int_{x_j}^{x_{j+1}} f(\varphi_{j}-B_{j+1})\ dx, \end{equation}.

Similar with  a special case of quadratic UPG presented in  \cite{connections4CD},   
we have: 
\begin{theorem}\label{thm:PGE-qBs}
Assume that   $f \in C([0,1])$  and $u_h=\sum_{j=1}^{n-1} u_j\, \varphi_j$ is  the solution of the UPG formulation \eqref{eq:1d-modelPG}, with  the bubble  function $B$  chosen such that  $b=\frac12 -\frac{\varepsilon}{h}>0$.  Then, 
\[
 \left | u_j - \int_0^{x_j} f(x)\ dx\right |  \leq 2\, \|f\|_{\infty}  \left (1 -\frac{\varepsilon}{h} \right ) h 
 \leq 2\, \|f\|_{\infty}\, h.
\]
\end{theorem}
\begin{proof}
We note  that \[\int_0^{x_1} B_1\,dx  =  \int_{x_j}^{x_{j+1}} B_{j+1}\, dx = \left (\frac12 -\frac{\varepsilon}{h} \right)\, h,  \ \text{and} \  \int_0^{x_1} \varphi_1\ dx  =   \int_{x_j}^{x_{j+1}}  \varphi_j=\frac{h}{2}.
\]
Thus,  using the formulas \eqref{eq:eplpicitSol} and the  triangle inequality, we have
\[
\begin{aligned}
 \left | u_j - \int_0^{x_j} f(x)\ dx\right | &= \\
 \left |  \int_0^{x_1} f(B_1-\varphi_{0})\ dx\,  + \int_{x_j}^{x_{j+1}} f(\varphi_{j}-B_{j+1})\ dx\right | & \leq \|f\|_{\infty}  \left (2 -\frac{2\, \varepsilon}{h} \right ) h.
\end{aligned}
\]
\end{proof}
Consequently, for $ j=1,\cdots, n-1$, the components $u_j$ of the UPG solution approximate $w(x_j)= \int_0^{x_j} f(x)\ dx$ with order $\mathcal{O}(h)$. 
If  $f$ is independent of $\varepsilon$ and $h >  \varepsilon$, then the UPG solution  is $\mathcal{O}(h)$ close to the solution of the transport problem \eqref{VF1d-reduced} at the interior nodes, hence  free of non-physical oscillations. As mentioned before,  for $\varepsilon \ll h$ and  for $ j=1,2, \cdots, n-1$, 
the values $w(x_j)$ are very close to $u(x_j)$, where   $u$ is the exact solution  of \eqref {eq:1d-model}. Thus, at least 
for   $\varepsilon \ll h$, the discretization approximates well the exact solution at the interior nodes. 
In Section \ref{sec:UPGinfError}, we show  that, for  different  generating bubbles $B$,  we can provide  precise  discrete infinity norm error estimates for the UPG  method.


\section{The importance of  the discrete infinity error control} \label{sec:UPGinfError}

Many finite element discretizations of convection dominated problems, including the Streamline Upwind Petrov–Galerkin (SUPG),  SPLS and Stream line Diffusion (SD), exhibit non-physical oscillations.  
For such methods, the convergenge analysis  was  done in various weighted energy norms that  include $L^2$ and  $H^1$  terms and the  error estimates  were  done on the  entire domain of the problem, including on regions corresponding to boundary layers. 
When considering such energy norms, the best approximation of the  finite element solution is usually estimated by  using the approximation properties  of the interpolant in $L^2$  and $H^1$ norms.  The interpolant approximation error  could be large,  especially  on subdomains containing the boundary layer. 


As an example of possible  large magnitude for  $|u-I_h(u)|$,  we considered $f=1$ in  \eqref{eq:1d-model}  with  the exact solution  
$
u(x)= x- (e^{\frac{x}{\varepsilon}} -1)/(e^{\frac{1}{\varepsilon}} -1). 
$
As presented in  \cite{UPG4CD}, the  linear  interpolant of $u$ on the uniform mesh of size $h$ satisfies 
\[
\begin{aligned}
{|u- I_h(u)|_{[0, 1]}^2} &=  \frac{1+e^{-1/\varepsilon}} {1-e^{-1/\varepsilon}} 
\left (\frac{1}{2 \varepsilon} - \frac{1}{h} \, \frac{1-e^{-h/\varepsilon}}{1+e^{-h/\varepsilon}} \right ) 
 {\approx \frac{1}{2\varepsilon} -\frac{1}{h}, \ \text{for} \ \varepsilon \ll h}, \\
|u- I_h(u)|_{[0, 1-h]}^2 & =  \frac{e^{-2h/\varepsilon}-e^{-2/\varepsilon}}{1-e^{-2/\varepsilon}} \, 
|u- I_h(u)|_{[0, 1]}^2  \approx e^{-2h/\varepsilon} |u- I_h(u)|_{[0, 1]}^2. 
\end{aligned}
\]
Thus, for $\varepsilon \ll h$, we have
\[
\begin{aligned}
|u- I_h(u)|_{[0, 1-h]} &  \approx e^{-h/\varepsilon} |u- I_h(u)|_{[0, 1]} \approx 0, \ \text{and}  \ \\
|u- I_h(u)|^2_{[1-h, 1]} & \approx |u- I_h(u)|^2_{[0, 1]} \approx \frac{1}{2\varepsilon} -\frac{1}{h}. 
\end{aligned}
\]
These estimates show that, for $ \varepsilon \ll h$,  the interpolant's energy error \\  $|u-I_h(u)|$ 
is insignificantly small on the  interval $[0, 1-h]$, is very large on $[0, 1]$, and is essentially attained  on the last mesh-interval $[1-h, 1]$. Thus, energy type error  estimates  are not suitable to provide higher order of approximation of the solution.  Therefore, they are not able to predict non-physical  oscillations of the discrete solution. 

 

 
 \subsection{The importance of convergence in discrete infinity error norm } \label{asec:ConvDiscInf} 
 The idea of  controlling the discrete infinity error in the convergence analysis  for convection diffusion problems can be found in the context of  monotone schemes discretization, as  presented for example in   \cite{XuZiEAFE}. 
 
 In order for a  continuous  piecewise linear  approximation  $u_h$  of   the solution $u$ of  a convection dominated problem to be oscillation free, it would be enough to have that {$\|u_h -I_h(u)\|_{h,\infty} $} converges fast to zero, as $h\to 0$. Here,
 \[
 \|u_h -I_h(u)\|_{h,\infty}:=\max\{|u_h(x_j)-u(x_j)| : \, x_j \in \Omega, \ x_j -\text{mesh node}\}.
 \] 
 This is also known as {\it the closeness property} of the discrete solution. 
 
 For example, let us  assume that for a particular discrete  solution \\  $u_h\in \M_h=C^0-P^1$, we have 
  \[
 {\|u_h -I_h(u)\|_{h,\infty} \leq \mathcal{O}(h^2)}.
\]
For functions in $\M_h$, by using standard inverse inequality  estimates that  hold with constants independent of $h$ and $\varepsilon$,  we also have 
\begin{equation}\label{eq:InvIneq}
 h\, |u_h -I_h(u)| \preceq \|u_h -I_h(u)\| \preceq  \|u_h -I_h(u)\|_{h,\infty}.   
\end{equation}
Consequently,  the splitting $u_h-u = (u_h- I_h(u)) + (I_h(u)-u)$ leads to {optimal error estimates} 
 for $|u_h-u|$ and $\|u-u_h\|$  {\it away from the boundary layers} (ABL). Here, we can  define ABL as any  subdomain region on which the   $H^1$ or $L^2$ interpolation  errors maintain optimal order  for the exact solution, and the error estimates are independent of $\varepsilon$.   
 
 In conclusion, controlling the discrete infinity error eliminates oscillations and could lead to optimal $L^2$ and  $H^1$ error estimates, ABL. The norm estimate \eqref{eq:InvIneq} is independent of the problem dimension, hence the ideas of this section can be extended  to {\it multidimensional} convection-diffusion problems. 



\subsection{The exponential bubble  UPG provides the exact solution at the nodes} 

In \cite{connections4CD},  we analyzed an UPG discretization for \eqref{eq:1d-model}, using the  special  exponential generating bubble function  $B:[0, h] \to \R$ defined by 
 \begin{equation}\label{expB}
  {B(x)=B^e(x):=\frac{1 - e^{-\frac{x}{\varepsilon}}}{1 - e^{-\frac{h}{\varepsilon}}} - \frac{x}{h}}.  
  \end{equation}
Introducing the notation  ${t_0:=\tanh\left (\frac{h}{2\varepsilon}\right )}$,  we get  $b= \frac{1}{h}\,  \int_0^{h} B(x)\, dx = \frac{1}{2 t_0} - \frac{\varepsilon}{h}$ and  
   \begin{equation}\label{IntexpB}
 M_{fe} = M^e_{fe} = \frac{1}{{t_0}} tridiag\left ( -\frac {1+{t_0}}{2} ,\ 1,\   -\frac{1-{t_0}}{2} \right).
     \end{equation}
It has been shown  that  the exponential bubble UPG recovers  the exact solution $u$ at the mesh nodes, see   \cite{connections4CD} and   \cite{roos-stynes-tobiska-96}. 

Next, we emphasize that, for $ \varepsilon \ll h$, at the interior nodes, the exponential bubble UPG solution is also very close  to the solution of the LR transport problem \eqref{VF1d-reduced}. 
 For $\frac{\varepsilon}{h} \to 0$, we have {\it with fast convergence} that 
\[ 
t_0 \to 1, \ \text{and} \  g_j =\varphi_j +B^e_j-B^e_{j+1} \to \chi_{|_{[x_{j-1}, x_j]}}.
\]
Consequently, we have  {\it with fast convergence} that  
\[
M_{fe} \to  tridiag(-1, 1, 0), \ \text{and} \ (f,g_j) \to  \int_{x_{j-1}}^{x_j} f(x) \, dx. 
\]
Thus, for $ \varepsilon \ll h$,  the matrix $M_{fe}^e$ is very close to  $tridiag(-1, 1, 0)$, 
and the {\it reduced  linear system}   becomes 
\[
 [tridiag(-1, 1, 0)] \, W = \left [ \int_{x_{0}}^{x_1} f(x) \, dx, \cdots,  \int_{x_{n-2}}^{x_{n-1}} f(x) \, dx \right ]^T.
\]
By forward solving the linear system, we obtain
\[
w_j = \int_{0}^{x_j} f(x) \, dx, \ j=1,2,\cdots,n-1. 
\]
This implies that the component $u_j$  of the UPG discrete solution is very close or identical to the value $w(x_j)= \int_{0}^{x_j} f(t) \, dt$. 
In conclusion, we have:
\begin{remark} \label{rem:RS}
The  {\it reduced  linear system} obtained from taking the limit as  ${\varepsilon}/{h} \to 0$ in the linear system for the exponential bubble UPG method,  represents  a discretization  of  the {\it  LR transport problem}  \eqref{VF1d-reduced}, which has unique solution.  For $ \varepsilon \ll h$, at the interior nodes, the exact solution  is very close to the solution  $w$  of \eqref{VF1d-reduced}. 
 \end{remark}
 
 Another benefit of the exponential bubble UPG method, 
is that,  using the Green's function for the  problem  \eqref{eq:1d-model}, i.e., 
  \[ 
 G(x,s)=\frac{1}{e^\frac{1}{\varepsilon}-1} \begin{cases}(e^\frac{1}{\varepsilon}-e^\frac{x}{\varepsilon})(1-e^{-\frac{s}{\varepsilon}}),&0\leq s < x,\\ (e^\frac{x}{\varepsilon}-1)(e^\frac{1-s}{\varepsilon}-1),&x\leq s\leq 1, \end{cases}
 \] 
we have  an {\it exact  formula for the  inverse of the UPG discretization matrix}: 
   \begin{equation} \label{eq:iMUPG}
{M_{fe}^{-1} = \left [ G(x_j, x_i) \right ]_{i,j=1,2,\cdots,n-1}}.
\end{equation}
Consequently, for any inside node $x_j$, the UPG solution $u_h$ is given by
 \begin{equation} \label{eq:expF}
 {u_h(x_j)=u(x_j)} = { \sum_{i=1}^{n-1} G(x_j, x_i) \left (f, \varphi_i + B^e_i-B^e_{i+1} \right)}.
\end{equation}
The formula allows for discrete infinity error analysis for other bubble UPG methods. Details of the proofs are given  in \cite{UPG4CD}. 

\subsection{Properly scaled quadratic bubble  UPG provides optimal approximation}
Working with exponential bubble functions could lead to difficulties  in computing the right hand side dual vectors.  A {\it properly scaled quadratic  bubble  UPG} can be used to obtain $\mathcal{O}(h^2)$  approximation of the exact solution in the  {\it discrete infinity norm}. 
We  choose the generating bubble function 
\[
 { B(x)= B^q(x):= \frac{4\, \beta}{h^2} x(h-x)}, \ \text{with} \   {\beta=\frac32 \left (\frac{1}{2 t_0}- \frac{\varepsilon}{h}\right )}, 
\]
such that the quadratic bubble  $B^q$ and the exponential bubble $B^e$ have the same average $b$. The stiffness matrix formula \eqref{eq:M4CDfe} gives that the  quadratic and  exponential UPG discretizations lead to the same linear system matrix. From \eqref{eq:iMUPG}, we get that for any inside node $x_j$, the quadratic bubble solution $u_h$ is given by
 \begin{equation} \label{eq:qF}
u_h(x_j)= { \sum_{i=1}^{n-1} G(x_j, x_i) \left (f, \varphi_i + B^q_i-B^q_{i+1} \right)}.
\end{equation}
By using the solution formulas \eqref{eq:expF} and \eqref{eq:qF} and standard linear algebra arguments, we obtain the following result. 
\begin{theorem}  If $\displaystyle u_h= \sum_{j=1}^{n-1} u_j \varphi_j$ is the {quadratic bubble} UPG solution with 
 $\beta=\frac32 \left (\frac{1}{2 t_0}- \frac{\varepsilon}{h}\right )$, and  {$e^{-\frac{\varepsilon}{h}} \leq h$}, then
\[
 {\|u_h -I_h(u)\|_{h,\infty}:=\max  |u(x_j)-u_j | \leq    2 \varepsilon \|f\|_\infty + \frac{h^2}{4} \|f'\|_\infty}.  
\]
\end{theorem}

The proof can be found in  \cite{UPG4CD}. 
 The estimate is a  sharp result, especially if  $\|f\|_\infty $ and   $\|f'\|_\infty$ can be bounded independently of $\varepsilon$. 
 
 For  $\varepsilon  \leq  h^2$, the  theorem gives  $ {\|u_h -I_h(u)\|_{h,\infty}  \leq C h^2}$, 
with precise control of the constant $C$. In addition, according to  Section \ref{asec:ConvDiscInf}, the estimate leads to 
\[
{ |u-u_h| \leq \mathcal{O}(h)} , \text{and} \ {\|u-u_h\|_{L^2}  \leq \mathcal{O}(h^2)}\ \text{{away} form the boundary layer}.
\]
\begin{remark}
We note that  it is essential to use the scaling  $\beta =\frac32 \left (\frac{1}{2 t_0}- \frac{\varepsilon}{h}\right )$ in order to obtain optimal approximation in the  discrete infinity norm. 
 Otherwise, the    $H^1$ and $L^2$ errors  increase and  the orders of convergence decrease,  especially for the $L^2$ error. 
\end{remark} 
The strategy of  properly scaling the generating bubble function could be used to obtain  {\it  discrete infinity error  estimates for  the two dimensional  case} as well.


\section{Two dimensional  extensions  of bubble UPG discretizations} \label{ss:PG2D}

The bubble UPG  idea can be extended to multidimensional case of problem \eqref{eq:2d-model}.  In \cite{UPG4CD}, we described the extension ideas for $\Omega=(0,1)\times (0, 1)$ and $\textbf{b}=[1,0]$.   In this case, the problem  \eqref{eq:2d-model} becomes: Find $u=u(x,y)$ such that 
 \begin{equation}\label{eq:b10}
 \left\{
\begin{array}{rcl}
     -{ \varepsilon}\, \Delta u +  u_x & =\ f & \mbox{in} \ \ \ \Omega,\\
      u & =\ 0 & \mbox{on} \ \partial\Omega.\\ 
\end{array} 
\right.
 \end{equation} 
It is known that the general solution of \eqref{eq:b10} could exhibit both elliptic and parabolic boundary layers, see e.g.,  \cite{roos-stynes-tobiska-96}. 
The following general  bubble  UPG in the $x$-direction discretization   was  considered in \cite{UPG4CD}.  With the notation of Section  \ref{sec:PG}, using a  generating  bubble  $B=B(x)$, we  define the trial and test spaces by: 
\[
\M_h= \text{span}\{\varphi_i(x)\,\varphi_j(y),  \ \text{for all } (x_i,x_j) \in  \Omega  \}, \  \text{and}
\] 
\[ 
 V_h= \text{span}\{(\varphi_i(x)+ B_i(x)-B_{i+1}(x)) \,\varphi_j(y),  \ \text{for all } (x_i,x_j) \in  \Omega \}. 
\]
A general UPG  discretization with  bubble functions in the $x$-direction for solving \eqref{eq:b10} is: Find $u_h \in \M_h$ such that 
\begin{equation}\label{eq:UPGxV}
b(v_h, u_h):= \varepsilon\, (\nabla u_h, \nabla v_h)+  \left (\frac{\partial u_h}{\partial x},v_h\right ) =(f,v_h) \ \Forall  v_h \in V_h. 
\end{equation}
 For defining the UPG with  {\it quadratic bubble upwinding}, we choose  
\[
 B(x)= B^q(x)= \frac{4\, \beta}{h^2} x(h-x), \, \text{with the matching} \   \beta=\frac32 \left (\frac{1}{2 t_0}- \frac{\varepsilon}{h}\right). 
\]
In this section,  we will focus on understanding the oscillation spikes  for the  quadratic UPG discretization of  \eqref{eq:b10} along the parabolic  boundary layers. 

Numerical tests  for examples with elliptic Boundary Layers (BL), show  that,  for $\varepsilon \leq h^2$ as in the one dimensional case,  the {\it UPG solution}  $u_h$,  satisfies 
\begin{equation}\label{eq:DinfErO2}
 {\|u_h -I_h(u)\|_{h,\infty} = \mathcal{O}(h^2)}, \ \text{and}
\end{equation}
\begin{equation}\label{eq:L2H1Er}
 |u-u_h| = \mathcal{O}(h) , \text{and} \  \|u-u_h\|_{L^2}  =\mathcal{O}(h^2), 
\end{equation}
away from the {\it  elliptic  boundary layer},  see \cite{UPG4CD}.
In addition, if the  $H^1$ and $L^2$ errors  are computed on the whole domain, we notice that 
the orders of convergence decrease. However,  in this case of  for $\varepsilon \leq h^2$, we get 
\begin{equation}\label{eq:L2H1Int}
{ |u-u_h| \approx  | I_h(u)-u_h| } , \text{and} \ {\|u-u_h\|_{L^2}  \approx \|I_h(u)-u_h\|_{L^2} }.
\end{equation}

Numerical tests for examples with elliptic BL near $x=1$  and parabolic BLs near  $y=0$ and $y=1$ showed that  the estimate \eqref{eq:DinfErO2} also  holds  for $\varepsilon \leq h^2$, but only {\it away from the parabolic BL}, and that  \eqref{eq:L2H1Er} still holds away from both elliptic and parabolic boundary layers. 

The numerical solution for $f=1$ (for which the solution exhibits both types of BL) has  {\it no spikes  for $\varepsilon >  h^2$}, but exhibits  {\it non-physical spikes} for $\varepsilon \ll h^2$ near  $y=0$ and $y=1$,  see Figure \ref{fig:f=1proj} (Left plot). 
Other discretization methods, including SUPG, show similar behavior for  $\varepsilon \ll h^2$.  In \cite{46Knobloch} and many of the references therein,  various stabilization methods were introduced  in order  to try  eliminating  such oscillations. 

\begin{figure}[h!]
	\begin{center}
		\includegraphics[height =1.8in, width=2.7in]{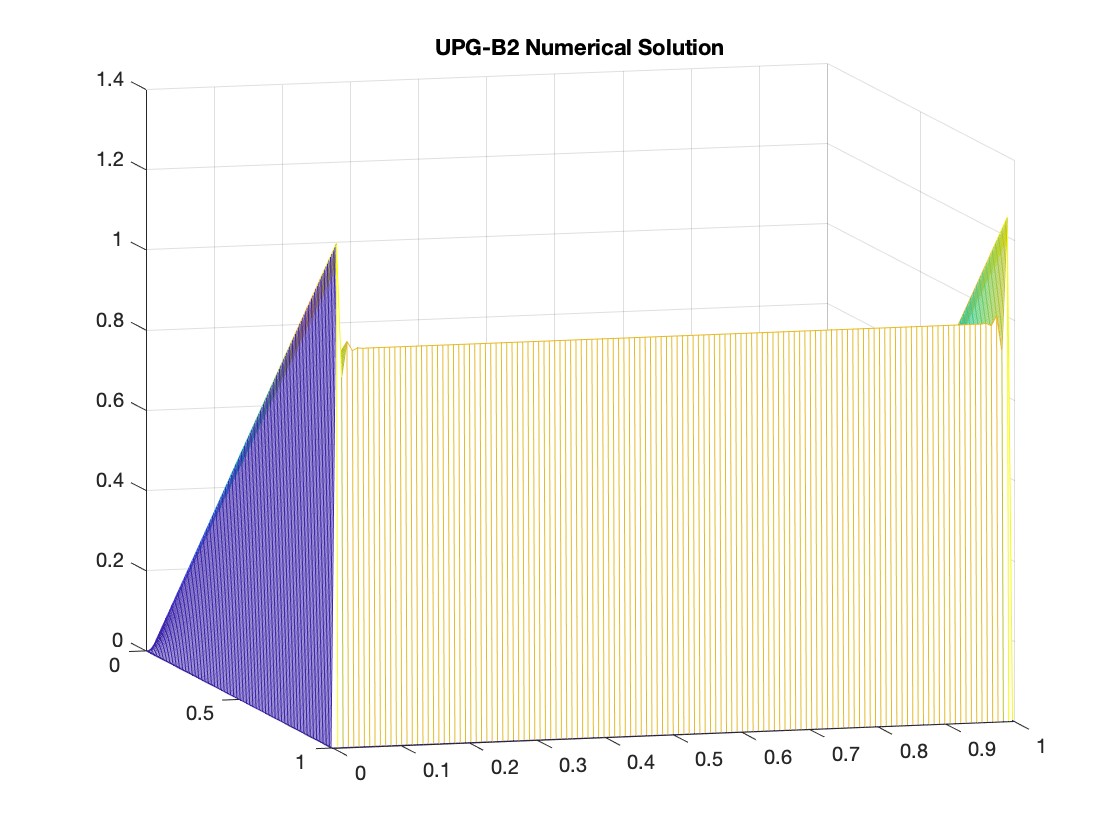}
		\includegraphics[height =1.2in, width=2.23in]{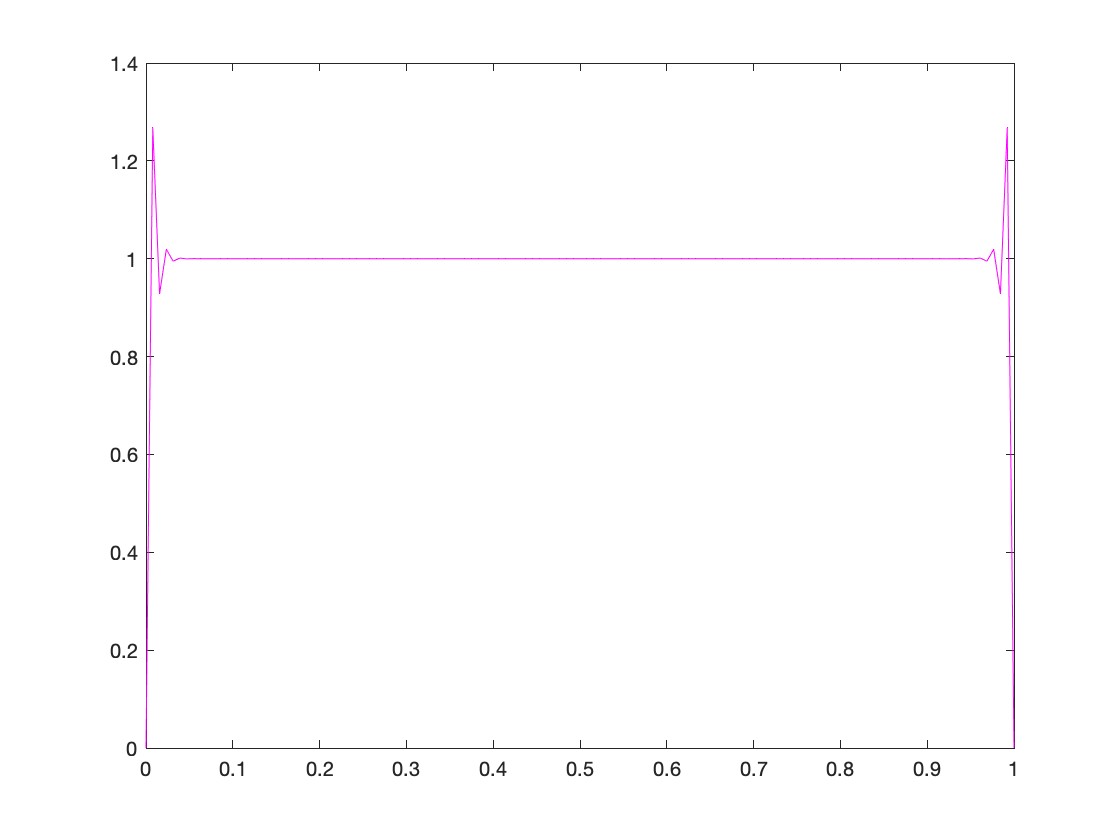}
	\end{center}
	\caption{The  quadratic UPG solution, $\varepsilon=10^{-7} \\ h=1/2^7$, and the $L^2$ orthogonal  projection of $1_{|_{[0, 1]}}$ on $\M_h$.}
		\label{fig:f=1proj}
\end{figure}	


\subsection{Understanding the oscillation spikes near the parabolic BLs}
We provide an  explanation of spikes along parabolic BLs for the 2D scaled {\it quadratic bubble UPG method}, based on  investigating the tensor structure of the resulting  linear system.  The linear system corresponding to  \eqref{eq:UPGxV} is 
\begin{equation}\label{eq:sys2dq}
A^q\, U^q = F^q, \ \text{where} \ A^q = M\, \otimes M^e_{fe} + {\frac{\varepsilon}{h}} \, S \otimes M^q. 
\end{equation}

Here,  $M=[(\varphi_i, \varphi_j)]_{i,j=1,\cdots, n-1} =  \frac{h}{6} \, tridiag(1, 4, 1)$ is the one dimensional mass matrix,
\[
 S=tridiag(-1, 2, -1),  \ \text{and} \  M^q = M+\beta\, ({h}/{3}) \, tridiag(-1, 0,1), 
 \]
\[
\begin{aligned}
U^q&=[u_{11}, \cdots, u_{n-1,1}, u_{12},u_{22}, \cdots,u_{n-1,2},\cdots, u_{n-1,n-1}]^T,\\
F^q&=[(f,g^q_1\,\varphi_1),\cdots, (f,g^q_{n-1}\, \varphi_1), (f,g^q_1\, \varphi_2), \cdots, (f,g^q_{n-1}
\, \varphi_{n-1})  ]^T, 
\end{aligned}
\]
$u_{ij}$ corresponds to the node $(x_i, x_j)$, and  $g^q_{i} =  \varphi_i + B^q_i -B^q_{i+1}$. 

Next, we note that  
\[
 \lim_{\varepsilon/h \to 0} M^e_{fe} = {C^0:=tridiag(-1, 1,0)},  \ \text{and} \lim_{\varepsilon/h \to 0} B^q \to B^{q_0}:=3\frac{x}{h} \left (1 -\frac{x}{h}\right).
 \]
By taking the limit as ${\varepsilon/h \to 0}$ in the system \eqref {eq:sys2dq},  we obtain the {\it reduced linear system}
\begin{equation}\label{eq:RLS}
\left [ M \, \otimes C^0\right ] W = F^{q_0}, \ \text{where} 
\end{equation}
\[
F^{q_0}=[(f,g^{q_0}_1\,\varphi_1),\cdots, (f,g^{q_0}_{n-1}\, \varphi_1), (f,g^{q_0}_1\, \varphi_2), \cdots, (f,g^{q_0}_{n-1}
\, \varphi_{n-1})  ]^T,
\]
and $g^{q_0}_{i} =  \varphi_i + B^{q_0}_i -B^{q_0}_{i+1}$.  

For  { $\varepsilon  \ll h$},  the solution $U^q$ of 
\eqref {eq:sys2dq} is very close to the solution $W$ of the {\it reduced linear system} \eqref{eq:RLS}. The solution of the linear system \eqref{eq:RLS} is 
\begin{equation}\label{eq:RLSi}
W= \left [ M^{-1} \, \otimes (C^0)^{-1} \right ]\, F^{q_0}. 
\end{equation}
By using that 
  \[
\int_0^{x_1} B_1^{q_0}(x)\, dx= \int_{x_j}^{x_{j+1}} B_{j+1}^{q_0}(x)\, dx  = \frac{h}{2}= \int_0^h \varphi_0(x)\, dx = \int_{x_j}^{x_{j+1}} \varphi_{j+1}\, dx, 
\]
and  \eqref{eq:RLSi}, for $f=1$ we obtain that, for  each $i=1, \cdots, n-1$ the $x_i$-section of  $W$  satisfies 
\[
  [w_{i1}, \cdots,  w_{i,n-1}]^T  = M^{-1}  \left[  \left (\int_0^{x_i} 1 \,ds,\varphi_1(y)\right ), \cdots,   \left(\int_0^{x_i}1 \,ds,\varphi_{n-1}(y)\right ) \right ]^T.
\]
 Since $M$ is the one dimensional mass matrix,  the  above inversion  represents  the projection of the  constant function  $x_i\, 1_{|_{[0, 1]}}$ to $span\{\varphi_1(y), \cdots, \varphi_{n-1}(y)\}$.  Figure \ref{fig:f=1proj}, shows the quadratic UPG solution for a case  $\varepsilon  \ll h$ versus the one dimensional  $L^2$ projection of
the constant function  $1_{|_{[0, 1]}}$ onto $\M_h$. It is interesting to see that the graph  of the $L^2$ projection  is scaled down to generate the quadratic UPG solution, when $x$  goes from $1-h$ to $0$. This shows   that  the UPG  oscillations in the  case  $\varepsilon  \ll h$ are due to the one dimensional projection of a function that is not zero at the ends onto the space $\M_h$ of continuous piecewise linear functions {\it with zero boundary conditions}. 

 We also note  that the the $x_i$-section of  $W$, the solution of the {\it reduced linear system} \eqref{eq:RLS}, is approximated by  the projection of $\int_0^{x_i} f(s,y)\, ds$ on $\M=\M(y)$. Thus, $W$ corresponds to a discretization of the {\it continuous transport problem}: Find $w \in H^1(\Omega)$ such that
\begin{equation}\label{eq:b10r}
 \left\{
\begin{array}{rcl}
      w_x & =\ f & \mbox{in} \ \ \ \Omega,\\
      w & =\ 0 & \mbox{for} \ x=0, y=0, \ \text{and}\  y=1.\\ 
\end{array} 
\right.
 \end{equation} 
If $f\in C^1([0 1] \times [0, 1]) $ and  $f(x,0) =f(x,1)=0$  for all $x \in [0, 1]$, then the exact solution  of  \eqref{eq:b10r} is $w(x,y)=\int_0^xf(s,y)\, ds$ for al $(x,y)\in \Omega$. On the other hand, if  $f\in C^1([0, 1] \times [0, 1]) $, but  $ f(x,0) \neq 0$ or $f(x,1) \neq 0$ for some $x \in [0, 1]$, for example $f=1$, then \eqref{eq:b10r} might not have a  solution in $H^1(\Omega)$. 




\section{Conclusions on oscillatory behavior of convection-diffusion discreretization}\label{sec:conclusion}  
This paper  addresses   the oscillatory behavior of certain  finite element discrezitaions for a model  convection-diffusion  problem. The work identifies  the  causes  of non-physical oscillations of  finite element approximation  of convection dominated problems, and suggests  ways  to avoid such behavior.  For  the bubble UPG method, we emphasize on a new approach for conducting error analysis that starts with establishing  the {\it closeness} of the discrete solution first, followed by  using inverse inequalities to establish approximation estimates  in standard $L^2$ and $H^1$ norms, away from boundary layers. 
The  ideas presented here can be used in building  new discretizations that are free of non-physical oscillations for the multi-dimensional  convection dominated problems.  
 \vspace{0.1in}
 
 Here are the  conclusions that follow from the work in this paper, intertwined with   related work presented  in  \cite{UPG4CD, CRD-results, connections4CD, comparison4CD}. 


\begin{enumerate}
 \item  [1)] The behavior of  {\it the reduced discrete problem}  obtained by letting \\ $\varepsilon \to 0$ in the discrete variational formulation, 
 {\it predicts}   {\it  the numerical pollution} for the case $\varepsilon  \ll h$.
  \vspace{0.1in}

\item [2)]  To avoid  oscillations for discretizations of a convection dominated problem, {\it the reduced discrete problem} 
should have  unique solution  and should correspond to a  continuous problem that has unique solution, for example, a related transport equation.  
  \vspace{0.1in}

\item  [3)]  In the one dimensional case, to eliminate the non-physical solutions for  the  {\it standard linear} discretization or the {\it saddle point least square} discretization, one can split the data $f=(f-\overline{f}) + \overline{f} $, and solve the two corresponding  problems for $f-\overline{f}$ and $\overline{f}$. The solution for a constant data $\overline{f}$ can be found explicitly.  The data $f-\overline{f}$ has average zero, and both continuous and discrete corresponding  reduced problems have unique solutions. 

  \vspace{0.1in}   
    
   \item  [4)]  Local oscillations at the boundary of the domain, can appear because
the discrete solution is approximated by an $L^2$- projection of continuous functions with non-zero boundary conditions, to conforming subspaces that do not account for the boundary conditions.
      \vspace{0.1in}
    
\item [5)] For  the one dimensional quadratic  bubble UPG method, the special scaling of the bubble  leads to an optimal convergence estimate in the discrete infinity norm, and leads to  optimal error estimates in the standard  $H^1$ and $L^2$  norms, away from the boundary layers.
    \vspace{0.1in}
    
\item [6)] For the  two dimensional bubble UPG and a given $\varepsilon$, by choosing $h$ such that  $ h^2 \approx \varepsilon$,  the discrete solution is free of non-physical oscillations and approximates well the exact solution in both   $H^1$ and $L^2$ norms, {\it away from all boundary layers}. 

    \vspace{0.1in}
    
   \item  [7)] {The discretization of multi-dimensional  convection dominated problems}, could {\it benefit from the efficient discretization of the corresponding one dimensional  problem along each stream line}. We can construct discrete test spaces by a ``tensor type'' construction using 
 {UPG discretization  along the stream line direction}  and   {\it standard discretizations on the orthogonal direction(s)}.  
   \vspace{0.1in}
     
\item   [8)]  To eliminate non-physical oscillations, it is  desirable to have a discretization that satisfies a {\it closeness  property} such as 
\[
{\|u_h -I_h(u)\|_{h,\infty} \approx \mathcal{O}(h^\alpha)},  \text{with} \ \alpha >1. 
\]
Once this is achieved,  the estimate can be combined with standard inverse inequalities to prove   error estimates for $ |u-u_h|$ and  \\ $\|u-u_h\|_{L^2}$,  {\it away from  the boundary layers}. 
  
\end{enumerate}



 \bibliographystyle{plain} 
 \def\cprime{$'$} \def\ocirc#1{\ifmmode\setbox0=\hbox{$#1$}\dimen0=\ht0
  \advance\dimen0 by1pt\rlap{\hbox to\wd0{\hss\raise\dimen0
  \hbox{\hskip.2em$\scriptscriptstyle\circ$}\hss}}#1\else {\accent"17 #1}\fi}
  \def\cprime{$'$} \def\ocirc#1{\ifmmode\setbox0=\hbox{$#1$}\dimen0=\ht0
  \advance\dimen0 by1pt\rlap{\hbox to\wd0{\hss\raise\dimen0
  \hbox{\hskip.2em$\scriptscriptstyle\circ$}\hss}}#1\else {\accent"17 #1}\fi}

 \end{document}